

\magnification=\magstephalf
\hsize=156truemm
\vsize=240truemm

\def\al{\alpha}
\def\bet{\beta}
\def\gam{\gamma}
\def\del{\delta}
\def\ep{\varepsilon}

\def\la{\lambda}
\def\om{\omega}
\def\alfa#1{\ifx#10\alpha\fi\ifx#11\beta\fi\ifx#12\gamma\fi%
\ifx#13\delta\fi\ifx#14\ep\fi}
\def\al#1{\ifx#1_\expandafter\alfa\else\alpha #1\fi}
\def\Bbb{\bf} \def\frak{\bf}
\def\P{\Bbb P} \def\Z{\Bbb Z} \def\C{\Bbb C}
\def\cite#1{[#1]}
\def\Tr(#1){\mathop{\rm Tr}(#1)}
\def\Sing(#1){\mathop{\rm Sing}(#1)}
\def\Spec{\mathop{\rm Spec}}
\def\Pic{\mathop{\rm Pic}}
\def\roep #1.{\medbreak\noindent{\sl #1\/}.\enspace}
\def\endroep{\par\medbreak}
\def\wtx{\widetilde{X}}
\def\wl#1{\overline{#1}}
\def\sier#1{{\cal O}_{#1}}
\def\ha#1{{#1\over2}}
\def\square{\rm QED} \def\rtimes{\times}
\def\qed{\hfill$\qquad\square$\par\medbreak}

\catcode`[=\active
\def\macaulay{%
\catcode`[=\active
\def[##1]##2{\sb{##1}\ifx##22\sp\fi%
\ifx##23\sp\fi\ifx##24\sp\fi\ifx##25\sp\fi\ifx##26\sp\fi##2}
\parindent=0pt}
\catcode`[=12

\newcount\secno \newcount\subsecno

\font\largebf=cmbx10 scaled \magstep2
\font\largemath=cmmib10 scaled \magstep2
\outer\def\beginsection#1\par{\bigskip\bigskip\subsecno=0
\advance\secno by1
\message{#1}\leftline{\largebf \the\secno.\enspace#1}\nobreak
\medskip\noindent}

\def\subsec{\medskip\advance\subsecno by1
\noindent{\bf(\the\secno.\the\subsecno)}\enspace}
\long\def\subproclaim#1. #2\par{\medbreak\advance\subsecno by1
\noindent{\bf(\the\secno.\the\subsecno)\enspace#1.}\enspace{\sl#2}\par
\ifdim\lastskip<\medskipamount \removelastskip\penalty55\medskip\fi}
\def\subex#1{\medskip \advance\subsecno by1
\noindent{\bf(\the\secno.\the\subsecno)}\enspace{\sl Example\/}:
$n=#1$.\enspace}

\long\def\comment#1\endcomment{}

\def\biblio{\bigskip\bigskip
\message{Bibliography}\leftline{\largebf Bibliography}\nobreak
\medskip\frenchspacing
\parindent=0pt
\def\item##1##2\par{\ifx##1\undefined\else\smallskip\noindent
\rlap{[##1]}\hskip30pt
\hangindent=30pt ##2\fi}}
\def\rfarn{Ar} 
\def\rfbeht{Be} 
\def\rfbi{Bi} 
\def\rfbar{Ba} 
\def\rfcoh{Co} 
\def\rfdr{D--R} 
\def\rfgeer{Ge} 
\def\rfhal{Ha} 
\def\rfhir{Hi} 
\def\rfhul{Hul} 
\def\rfhur{Hu 1} 
\def\rfhurf{Hu 2} 
\def\rfkar{Ka} 
\def\rfklel{Kl} 
\def\rfkf{K--F} 
\def\rfmer{M\'e} 
\def\rfpel{Pe} 
\def\rfsb{SB} 
\def\rfjs{St} 
\def\rfvel{V\'e} 
\def\rfwei{We} 
\ifx\input amssym.def \input amssym\undefined\else \input amssym.def \input amssym \fi

\centerline{\largebf Degenerations of elliptic curves and cusp singularities}
\bigskip
\centerline{\bf Jan Stevens}
\bigskip
\bigskip

\begingroup
\narrower\narrower
\noindent{\bf Abstract.}\enspace
This paper gives more or less explicit
equations for all two-dimensional cusp singularities of
embedding dimension at least 4. They are closely related to
Felix Klein's equations for universal curves with  level $n$ structure.
The main technical result is a description of the versal deformation of an
$n$-gon in $\P^{n-1}$.
The final section contains the equations for smoothings of
simple elliptic singularities (of multiplicity $n\leq 9$).

\endgroup
\bigskip

Equations of cusp singularities are known for multiplicity at most $5$
\cite\rfkar.
Cusps are quite well understood in terms of generators
of the local ring, described as Fourier series;
the cusp is constructed as compactification
of the quotient of the product of upper half spaces
${\frak H}^2$ by the semi-direct product of a
lattice  in a real quadratic field and an infinite cyclic
group. The precise relations between these generators however are not known
even when equations for the cusp are available. Still the combinatorics of the
associated continued fractions provide the pattern for the equations.

The relation between cusp equations and families of elliptic curves
can already be seen in the hypersurface case. The cusps  of
multiplicity 3 are the $T_{pqr}$-singularities with $3\leq p\leq q \leq
r$ and at least one strict inequality. With a slightly different
notation I write
$$
X_0^{a_0+1}+X_1^{a_1+1}+X_2^{a_2+1}-X_0X_1X_2=0\;;
$$
the periodic continued fraction $[[a_0,a_1,a_2]]$ determines the cusp.
The equation  resembles the Hesse normal form
$$
X_0^3+X_1^3+X_2^3-3\mu X_0X_1X_2=0\;,
$$
which describes the universal elliptic curve with level $3$ structure.
For four values of the parameter ($\mu^3=1$ and $\mu=\infty$) the
curve degenerates into a triangle.
The projectivised tangent cone of the cusp is  the
triangle corresponding to $\mu=\infty$ in the Hesse normal form, if
all $a_i$ are at least 3. The other possibilities,
a line and a conic ($X_0^3-X_0X_1X_2$)
or a nodal cubic ($X_0^3+X_1^3-X_0X_1X_2$), do not occur as degenerate
curve with level 3 structure. Therefore it is better to forget about
extra structures and to directly consider
the versal deformation of the triangle:
$$
X_0X_1X_2-t_0X_0^3-t_1X_1^3-t_2X_2^3=0\;,
$$
in which the vertices can be smoothed separately.
The regular sequence $(t_0-X_0^{a_0-2},t_1-X_1^{a_1-2},t_2-X_2^{a_2-2})$
cuts out the cusp singularity from the total space of this
deformation.

For higher multiplicities I also start from equations
for elliptic curves with a level
structure. For odd $n$ they were written down by Felix Klein (see \cite
\rfkf).
The number of equations becomes rather large, but due to the symmetries of
the elliptic curve and the total family one needs only a few different
types. In fact, they can all be combined in the formula
$$
F_{\al_0\al_1\al_2\al_3}=
s_{\al_0-\al_1}s_{\al_2-\al_3}X_{\al_0+\al_1}X_{\al_2+\al_3} +
s_{\al_0-\al_2}s_{\al_3-\al_1}X_{\al_0+\al_2}X_{\al_3+\al_1} +
s_{\al_0-\al_3}s_{\al_1-\al_2}X_{\al_0+\al_3}X_{\al_1+\al_2}\;,
$$
with the $s_\al$  depending on the modulus of the curve.
For special
values of the modular parameter (at the `cusps') the elliptic curve
degenerates to an $n$-gon. I adapt the equations to give the versal
deformation of the standard $n$-gon (with the  as vertices
coordinate points in their natural order).
Equations of cusps are again
obtained by the specialisation $t_i=X_i^{a_i-2}$.
The change in the formulas for the modular curve involve the deformation
parameters $t_i$ making the coefficients $s_\al$ into power
series in the $t_i$. Therefore the resulting cusp equations are not longer
polynomial if $n>5$. For even $n$ similar results can be obtained, basically
with the same formulas.

For the special case of {\sl degenerate cusps\/}
Ruud Pellikaan found the equations  using
that the normalisation of a degenerate cusp consists of
cyclic quotient singularities \cite{\rfpel}.
A degenerate cusp is a non-normal surface singularity  which can also
be described by numbers $[[a_1,\dots, a_n]]$, but now
$a_i=\infty$ for at least one $i$; one formally puts $x_i^{\infty}=0$.
To describe the equations I use the abbreviation $X_i^{a_i-2}=t_i$ with
the convention that $t_i=0$ if  $a_i=\infty$.
\proclaim Proposition {\rm(Pellikaan)}.
Let $X(a_1,\dots,a_n)$ be a degenerate cusp with $n>3$. The\/ $n(n-3)/2$
equations
$$
F_{ij}\;\colon\quad
X_iX_j=X_{i+1}\left(\prod_{k=i+1}^{j-1}t_k\right)X_{j-1} +
       X_{j+1}\left(\prod_{k=j+1}^{i-1}t_k\right)X_{i-1}
$$
with $i-j\neq -1,0,1$ generate the ideal of the cusp.

Pellikaan's formulas applied
to isolated cusps with multiplicity $4$ and $5$ give the usual form
of the equations, but for higher multiplicities
one needs correction terms, which were already computed
by Pellikaan in the case $n=6$.
My cusp equations have a similar structure; they specialise to
both cases mentioned, of degenerate cusps and  elliptic curves.

The infinitesimal deformations of cusp singularities were determined
by Kurt Behnke using the Fourier series approach \cite\rfbeht.
In terms of my equations I can only give them for multiplicity at most
6; in those cases I  know also the versal deformation. I give the
formulas for $n=6$. Even for simple elliptic singularities I do not have
a general formula for infinitesimal deformations, but I did compute them
and also the
versal deformation in all cases where the singularity is smoothable
(i.e., up to $n=9$). A geometric description of the versal
deformation of simple elliptic singularities can be found in \cite\rfmer.

In the first Section I review the equations for elliptic curves and
their symmetries and indicate the modifications for the  case of
even degree. The examples introduce the equations which are used for
the computations of the versal deformation in the last Section. The
second Section contains the main technical result, the versal
deformation of an $n$-gon. The third Section applies it
to the cusp equations.

\beginsection Elliptic curves with level {\largemath n\/} structure

Equations for elliptic curves of degree $n$ were first written down
by Bianchi \cite{\rfbi} for small odd $n$ and for general odd $n$ by
Felix Klein \cite{\rfklel}. Their formulas describe the universal
family of elliptic curves with level $n$ structure (in fact as scheme
over $\Spec \Z[1/n]$ \cite{\rfvel}).  Hurwitz \cite{\rfhur}
realised that in the case of even $n$ basically the same embedding
works, although it does not give the  universal family.
The basic reference for these results is \cite{\rfkf, F\"unfter
Abschnitt};  see also \cite{\rfhul}.

\subsec
Let $E=\C/\Z\omega_1+\Z\omega_2$ be an elliptic curve and let
$n\geq5$ be an odd integer. Klein defined functions $X_\al(u)$,
$\al\in\Z/n$, whose zeroes  in the period parallelogram are
the points ${\al\over n}\omega_1+{i\over n}\omega_2$,
$i=0,\dots,n-1$. The values of these functions at the origin define
constants $s_\al:=X_\al(0)$.

\proclaim Proposition {\rm (Klein) \cite{\rfkf, V.1~\S9}}.
The functions $X_\al$ embed $E$ in $\P^{n-1}$. The ideal of $E$ is
generated by the quadratic equations:
$$
F_{\al_0\al_1\al_2\al_3}=
s_{\al_0-\al_1}s_{\al_2-\al_3}X_{\al_0+\al_1}X_{\al_2+\al_3} +
s_{\al_0-\al_2}s_{\al_3-\al_1}X_{\al_0+\al_2}X_{\al_3+\al_1} +
s_{\al_0-\al_3}s_{\al_1-\al_2}X_{\al_0+\al_3}X_{\al_1+\al_2}\;,
$$
with the indices  taking values in $\Z/n$.

\noindent
These equations are a consequence of the $\sigma$-relation:
$$
\displaylines{\qquad
\sigma(t+u)\sigma(t-u)\sigma(v+w)\sigma(v-w)+
\sigma(t+v)\sigma(t-v)\sigma(w+u)\sigma(w-u)\hfill\cr\hfill{}+
\sigma(t+w)\sigma(t-w)\sigma(u+v)\sigma(u-v)=0\qquad (1.)
\cr}
$$
about which  Weierstra\ss\ wrote: `Mann kann, ohne
von der Function $\sigma(u)$ irgend etwas zu wissen, direct
nachweisen, dass es eine vier willk\"urliche Constanten
enthaltende (transcendente) ganze Function der Ver\"anderlichen
$u$ giebt, welche f\"ur $\sigma(u)$ in die Gleichung $(1.)$
eingesetzt, dieselbe befriedigt' \cite{\rfwei} (for a
discussion of this functional equation see  \cite{\rfhurf} or
\cite{\rfhal, I p.~187}).

Each function $X_\al$ is defined as translate of the function
$\sigma(u;\omega_1,\omega_2/n)$ multiplied by a suitable
exponential factor, which I describe below for completeness; the equation
$F_{\al_0\al_1\al_2\al_3}$ is obtained from  $(1.)$ by the
specialisation $t={1\over2}u+{\al_0\over n}\omega_1$, \dots,
$v={1\over2}u+{\al_3\over n}\omega_1$. Each term requires the same
exponential factor. To define the $X_\al$ one starts with the
functions:
$$
\sigma_{\la,\mu}(u;\om_1,\om_2):=
e^{(\la \eta_1+\mu \eta_2)(u-\ha{\la \om_1+\mu \om_2})}
\sigma(u-\la\om_1-\mu\om_2;\om_1,\om_2)\;,
$$
where the constants $\eta_i$ are the
periods of the integral of the second kind. They  satisfy the
Legendre relation $\omega_1\eta_2-\omega_2\eta_1=2\pi i$ (I follow
Klein and Fricke and take the imaginary part of $\om_1/\om_2$
positive). If $(\wl \om_1, \wl \om _2)= (\om_1,\om _2/n)$, and
$\wl\eta_1$, $\wl\eta_2$ are the corresponding periods, then
$$
{\wl\eta_1\over \wl\om_1}-{n\eta_1\over \om_1} =
{\wl\eta_2\over \wl\om_2}-{n\eta_2\over \om_2} =: G_1
$$
and $X_\al$ is now defined as
$$
X_\al(u)=\rho\; (-1)^\al e^{-G_1u^2}
\sigma_{{\al\over n},0}(u;\om_1,{\om_2\over n})
$$
with $\rho$ a final term independent of $\al$ to make the $X_\al$
themselves modular forms with simple transformation behaviour.
Note that the functions $X_\al$ are defined for all
$\al\in\Z$, and satisfy  $X_{\al+n}(u)=X_\al(u)$.

\subsec
For even $n$ similar results are obtained with functions $X_\al(u)$
whose zeroes are not $n$-torsion points themselves, but shifted by
$u_0={1\over2}\omega_1+{1\over 2n}\omega_2$. More precisely,
$$
X_\al(u)=
\rho\; e^{-{\al\pi i\over 2n}-{5\pi i\over4}} e^{-G_1u^2}
\sigma_{{\al\over n},\ha1}(u;\om_1,{\om_2\over n})\;.
$$
The same
$\sigma$-relation (1.) now gives equations involving constants
$\sigma_{{\al\over n},0}(0;\om_1,{\om_2\over n})$ \cite{\rfkf,
p.~268}; a small computation shows that they are
proportional to $e^{\al \pi i/n}X_\al(u_0)$.
Therefore I define
$$
s_\al:=e^{\al \pi i/n}X_\al(u_0)\;.
$$

\proclaim Proposition. 
For even $n$ the ideal of $E$ in $\P^{n-1}$ is generated
by the functions $F_{\al_0\al_1\al_2\al_3}$, defined by the same
formula as for odd $n$, but with indices all in $\Z/n$ or
in ${1\over2}+{\Bbb Z}/n$.


\subsec
Different quartuples $(\al,\bet,\gam,\del)$ lead to the same equation.
To analyse the situation I use different indices and write
$$
F^{h}_{ijk}=
s_{k+j}s_{k-j}X_{h+i}X_{h-i} -
s_{k+i}s_{k-i}X_{h+j}X_{h-j} +
s_{j+i}s_{j-i}X_{h+k}X_{h-k}\;.
$$
The choice $h+i=\al_0+\al_1$, etc., leads to the following
equations (modulo $n$) determining
the transition from one set of indices to another:
$$
\eqalign{
\al_0+\al_1+\al_2+\al_3&=2h\cr
\al_0+\al_1-\al_2-\al_3&=2i\cr
\al_0-\al_1+\al_2-\al_3&=2j\cr
\al_0-\al_1-\al_2+\al_3&=2k\cr
}\qquad
\eqalign{
h+i+j+k&=2\al_0\cr
h+i-j-k&=2\al_1\cr
h-i+j-k&=2\al_2\cr
h-i-j+k&=2\al_3\;.\cr
}
$$
Putting $h-i=\al_0+\al_1$, etc., gives
the solution $(h-\al_0,h-\al_1,h-\al_2,h-\al_3)$ instead of
$(\al_0,\al_1,\al_2,\al_3)$. Other choices just lead to a permutation of the
indices.

For odd $n$ division by two with integer result is always possible, so it
suffices to take $h\in \Z/n$ and $0\leq i<j<k< {n\over2}$.
For even $n$ there are two types of equation, with $2h$ even or odd. In the
latter case all indices $h$, $i$, $j$ and $k$ are in ${1\over2}+{\Bbb Z}$.
As $(h,i,j,k)$ and $(h+{n\over2},{n\over2}-i,{n\over2}-j,{n\over2}-k)$ give
the same equation I can take $0\leq h<{n\over2}$ and $0\leq i<j<k\leq
{n\over2}$.

\subsec
\def\labelell{(\the\secno.\the\subsecno)}
The system of equations $F_{\al\bet\gam\del}$ admits a large symmetry
group.
The action of the  involution of the elliptic curve and translation
by points of order $n$ is given by the following formulas:
\medskip
\item{i)}
    $X_\al(-u)=(-1)^nX_{-\alpha}(u)$
\item{ii)}
    $X_\al(u+{\om_1\over n})=
        (-1)^ne^{\eta_1(u+{\om_1\over 2n})}X_{\al-1}(u)$
\item{iii)}
    $X_\al(u+{\om_2\over n})=
       (-1)^ne^{\eta_2(u+{\om_2\over 2n})}\ep^{-\alpha} X_\al(u)$,
\medskip
\noindent where $\ep=e^{2\pi i/n}$.
These formulas describe the action of the Heisenberg group, see
\cite{\rfhul, I.2.4}
The exponential factor is the same for all $X_\al$,
so in the action on $\P^{n-1}$
it drops out.

\subsec
To investigate the symmetries coming from the action of the
modular group the dependence of
the $X_\al$ on $\tau=\om_1/\om_2$ has to be considered.
For odd $n$ the
functions $s_\al(\tau):=X_\al(0;\omega_1,\omega_2)$ are  modular forms
for the principal congruence subgroup  $\Gamma(n)$  \cite{\rfkf, p.~280} and
they embed the modular curve $X(n)$ into ${\Bbb P}^{n-3\over2}$
(there are $n-1\over 2$ essentially different non zero $s_\al$, as
$s_\al=-s_{-\alpha}=-s_{n-\alpha}$ by formula i) in \labelell)
and the equations
$F_{\al_0\al_1\al_2\al_3}$ describe the universal elliptic curve with level
$n$ structure \cite{\rfbar}. For even $n$ one has $s_{\al+n}=-s_\al$, so in
this case $s_{n-\alpha}=s_\al$. The
number of essentially different constants is $n/2$.

In either case relations between the $s_\al$ are obtained by further
specialising the equations $F_{\al_0\al_1\al_2\al_3}$ or
$F^{h}_{ijk}$,
by putting $u=0$ for $n$ odd and $u=u_0$ for even $n$; in the latter case
one has to introduce an extra factor $(-1)$ if the index sum is bigger
than $n$: $X_iX_j$ specialises to $s_is_j$ if $0\leq i+j<n$, but to
$-s_is_j$ for $n\leq i+j<2n$.

{\def\al_#1{\ifcase#1 i\or j\or k\or l\fi}
\proclaim Lemma.
For $n\geq6$ the following equations hold:
$$
s_{\al_1-\al_0}s_{\al_3-\al_2}s_{\al_0+\al_1}s_{\al_2+\al_3} -
s_{\al_2-\al_0}s_{\al_3-\al_1}s_{\al_0+\al_2}s_{\al_3+\al_1} +
s_{\al_3-\al_0}s_{\al_2-\al_1}s_{\al_0+\al_3}s_{\al_1+\al_2}=0\;,
$$
with $0\leq\al_0<\al_1<\al_2<\al_3\leq n/2$.

If} $n=p$ is prime  these equations define the modular curve
$X(p)\subset \P^{n-3\over2}$ \cite{\rfvel},  but  they
do not necessarily generate the homogeneous ideal of the curve (see
the example $n=11$ below). For general odd $n$ the curve $X(n)$ is only
an irreducible component of the zero locus.

The action of the modular substitutions
$S=\left({1\atop0}{1\atop1}\right)$ and
$T=\left({0\atop1}{-1\atop0}\right)$ is \cite{\rfkf, V.2~\S 7}:
$$
\openup 2\jot
\def\al{\alpha}
\eqalignno{
S:\quad X_\al' &= \ep ^{-\ha{\al(n-\al)}}X_\al \cr
T:\quad X_\al' &= {i^\ha{n-1}\over \sqrt n}
\sum _{\bet=0}^{n-1}\ep^{-\al\bet}X_\bet\;.\cr
\noalign{\noindent For even $n$ one has \cite{\rfkf, V.2~\S 8}:}
S:\quad X_\al' &= \ep ^{{n\over 8}+\ha{\al^2}}X_\al \cr
T:\quad X_\al' &= {-1\over \sqrt n}
\sum_{\bet=0}^{n-1}\ep^{-\al\bet}X_\bet\;. \cr
}
$$
The $X_\al$ are now modular forms of level $2n$ or $4n$; if $n$ is
divisible by $4$, then the appropriate congruence subgroup is
$$
\Gamma(n,2n)=
\{\, \pmatrix{a&b\cr c&d\cr}\in SL(2,\Z) \mid
  a\equiv d\equiv 1 \pmod n,\quad  b\equiv c\equiv 0 \pmod{2n}\,\}\;.
$$
For $n$ the double of an odd integer one has the additional conditions
$$
a\equiv d\equiv 1 +k n,\quad b+ c\equiv nk(k+1) \pmod{4n}
$$
for the $X_\al$ to remain unchanged \cite{\rfkf, p.~289}, but they
only exclude substitutions which multiply all $X_\al$ by $-1$.

\proclaim Proposition.
The elliptic curves $\C/\Z\om_1+\Z\om_2$ and
$\C/\Z\om_1'+\Z\om_2'$ have the same image in $\P^{n-1}$ if and only if
$\om_1/\om_2$ and $\om_1'/\om_2'$ are equivalent under $\Gamma(n)$.

\roep Proof.
The case $n$ odd is contained in the results of Bartsch \cite{\rfbar}
and V\'elu \cite\rfvel, so
it suffices to look at the action of representatives of
$\Gamma(n)/\Gamma(n,2n)$ for even $n$. Up to a common factor the substitution
$S^n= \left({1\atop0}{n\atop1}\right)$ transforms the $X_\al$ into $(-1)^\al
X_\al$, whereas $u_0'=u_0+\ha n\om_2$, so also $s_\al$ is transformed
into $(-1)^\al s_\al$. The equations $F^h_{ijk}$ are
invariant: if all indices are integers, then
$h-i\equiv h+i\pmod 2$, whereas
$h-i\not\equiv h+i\pmod 2$ for elements of $1/2+\Z$.

One computes that $\left({1\atop n}{0\atop1}\right)=TS^{-n}T^{-1}$
sends (up to a common factor) $X_\al$ to $X_{\al-\ha n}$, whereas
$u_0'= u_0+\ha1\om_1$. Therefore the ideal is again invariant
and $s_\al'=s_\al$ as $X_\al'(u_0')=X_{\al-\ha n}(u_0+\ha1\om_1)=
X_\al(u_0)$ (up to a common factor).  \qed

\noindent The double cover of the modular curve $X(n)$
given by the involution $s_\al\mapsto (-1)^\al s_\al$ satisfies the
equations of the Lemma. So the equations $F_{\al\bet\gam\del}$ describe a
family of curves of genus one over $X(n)$.

\subsec
The equations are the $4\times 4$ Pfaffians of the anti-symmetric matrix
$P$ with entries $P_{\al\bet}=s_{\bet-\al{}}X_{\bet+\al{}}$. For odd $n$ one
obtains all equations in this way, while for even $n$ a second matrix is
necessary, say with $P_{\al\bet}=s_{\bet-\al{}}X_{\bet+\al+1}$, to yield
the equations with half-integer coefficients.
This way of writing the equations entails relations with linear
coefficients.

As the equations are not linearly independent (if $n>5$) one has also
relations whose coefficients depend only on the parameters $s_i$; they
involve the equations between them. It is not difficult to
write down a minimal set of generators of the ideal of the elliptic
curve, but then the  linear relations are much more complicated.

\proclaim Proposition.
The syzygies between the equations $F_{\al_0\al_1\al_2\al_3}$ of
the elliptic curve $E$ are generated by the relations:
$$
R^{h;i}_{jkl}:\quad
-s_{k+i}s_{k-i}F^{h}_{ijl}+
s_{l+i}s_{l-i}F^{h}_{ijk}+
s_{j+i}s_{j-i}F^{h}_{ikl}=0\;.
$$
and the linear relations:
$$
R_{\al_0;\al_1\al_2\al_3\al_4}:\quad
s_{\al_0-\al_1}X_{\al_0+\al_1}F_{\al_0\al_2\al_3\al_4} +
s_{\al_0-\al_2}X_{\al_0+\al_2}F_{\al_1\al_0\al_3\al_4}
+s_{\al_0-\al_3}X_{\al_0+\al_3}F_{\al_1\al_2\al_0\al_4}
+s_{\al_0-\al_4}X_{\al_0+\al_4}F_{\al_1\al_2\al_3\al_0}
=0\;.
$$

\noindent
The proposition follows from the analogous statement for
the more general equations of the degeneration which is the subject of the
next section.


\subsec
The curve $\frak H/\Gamma(n)$ can be compactified with
${1\over2}n^2\prod_{p|n}(1-{1\over p^2})$ cusps. Over a cusp
the  universal elliptic curve (for odd $n$) degenerates to an
$n$-gon. Each subgroup $\Z/n\subset \Z/n\times\Z/n$ has $n$
invariant hyperplanes
forming an $n$-simplex whose vertices are the fixed points of
the action on $\P^{n-1}$. For the subgroup generated
by  $(1,0)$ it is the simplex of reference,
with vertices $e_\al$.
By joining the vertices $e_\al$ and $e_{\al+m}$
by lines I obtain  an $n$-gon if $(m,p)=1$.
This construction yields $\phi(n)/2$ different $n$-gons.
Each
of them is a generalised elliptic curve whose origin is the
intersection with the $(n-3)/2$-dimensional subspace given  by
the equations
$X_\al+X_{n-\alpha}=0$; the intersection point lies on
the line joining $e_{n-m\over2}$ and $e_{n+m\over2}$ (consider the
indices as elements of $\Z/n$ or take  $m$ odd). Defining equations
for the $n$-gon are
$$
F^{h}_{0j{n-m\over2}}=s_{n-m\over2}^2X_{h+j}X_{h-j}\;.
$$
In the same way $n$-simplices can be constructed for the other of the
$n\prod_{p|n}(1+{1\over p})$ cyclic subgroups of order $n$
in $\Z/n\times\Z/n$; they are permuted by the modular group.

If $(n,m)=d$ then
the lines joining $e_\al$ to $e_{\al+m}$ in the standard simplex form
$d$ $n/d$-gons. The equations still reduce to
$F^{h}_{0j{n-m\over2}}=s_{n-m\over2}^2X_{h-j}X_{h+j}$; if $n/d=3$
they define $d$ planes (a plane triangle does not admit a quadratic
equation), otherwise the zero locus consists  $d$ $n/d$-gons. It
follows that the equations between the $s_\al$ have solutions other
than the curve $X(n)$, if $n$ is not a prime.

For an example of a degeneration with
$n$ even see the case $n=8$ below.

\comment
\medskip \noindent{\sl Example\/}: $n=5$.
There are five equations:
$$
\displaylines{
s_3s_1X_0^2-s_3^2X_1X_4+s_1^2X_2X_3 \cr
s_3s_1X_1^2-s_3^2X_2X_0+s_1^2X_3X_4 \cr
s_3s_1X_2^2-s_3^2X_3X_1+s_1^2X_4X_0 \cr
s_3s_1X_3^2-s_3^2X_4X_2+s_1^2X_0X_1 \cr
s_3s_1X_4^2-s_3^2X_0X_3+s_1^2X_1X_2 \cr}
$$
Eliminating the variables $s_1$ and $s_2$ gives a determinantal
surface $S_{15}$ of degree 15 whose normalisation is Shioda's modular
surface $S(5)$. For more details see \cite{\rfhul}.
\endcomment

\subex7
I write the equations with the $\sigma$-constants $s_1$, $s_2$ and
$s_4$; for every prime $p=4k-1$ it is convenient to use the
quadratic residues as indices. The modular curve is the Klein
quartic
$$
\macaulay
s[1]s[4]^3+s[2]s[1]^3+s[4]s[2]^3\;.
$$
The elliptic curves are defined by $14$ linearly independent
equations, out of $28$  quadrics $F^{h}_{ijk}$, of which I write
only the ones with $h=0$; the others can be
obtained by cyclic permutation of the $X_\al$.
$$
\macaulay
\catcode`*=9
\displaylines{
s[1]s[4]X[0]^2+s[2]^2X[1]X[6]-s[1]^2X[2]X[5]
\cr
s[1]*s[2]X[0]^2*+*s[4]^2X[2]*X[5]-*s[2]^2X[3]*X[4]
\cr
s[2]*s[4]X[0]^2+s[1]^2X[3]*X[4]-s[4]^2X[1]*X[6]
\cr
s[1]*s[4]*X[3]*X[4]+s[2]*s[4]*X[2]*X[5]+s[1]*s[2]*X[1]*X[6]
\cr}
$$
The modular substitution $U$ with $\om_1'\equiv 4\om_1$, $\om_2'
\equiv \om_2/4 \pmod{7}$  operates by
$X_\al'=X_{4\al{}}$ (see \cite{\rfkf, p.~302}), so it acts by
cyclic permutation of the indices $(1,2,4)$ and $(6,5,3)$; it
permutes the first three equations, whereas the fourth is invariant.
This reduces the number of essentially different equations to
two. Although all equations are  specified by only one formula
$F_{\al\bet\gam\del}$, there is no notable action of the available
symmetry groups which permutes the remaining two. However under the
action of $T$ one single  equation is transformed into a linear
combination of all $28$ equations.

\subex9
Again I write only equations of the form $F^0_{ijk}$. An appropriate
modular
substitution $U$  now permutes the indices $(1,4,7)$ and their complements
$(8,5,2)$ leaving $(0,3,9)$ invariant. The ten equations  split into four
types:
$$
\displaylines{
s_1^2X_2X_7-s_7^2X_1X_8+s_1s_3X_0^2 \cr
s_1^2X_3X_6-s_3^2X_1X_8-s_4s_7X_0^2 \cr
s_1s_3X_3X_6+s_4s_7X_2X_7-s_1s_4X_1X_8 \cr
s_1s_3X_4X_5+s_3s_4X_2X_7+s_3s_7X_1X_8\;. \cr}
$$
Notice that the last equation ($F^0_{124}$) can be divided by $s_3$,
to give  $s_1X_4X_5+s_4X_2X_7+s_7X_1X_8$.
As equations between the coefficients I get:
$$
\displaylines{
s_3^3s_1+s_7^3s_4-s_1^3s_4  \cr
s_3^3s_4+s_1^3s_7-s_4^3s_7  \cr
s_3^3s_7+s_4^3s_1-s_7^3s_1  \cr
(s_1^2s_7+s_4^2s_1+s_7^2s_4)s_3\;.  \cr}
$$
In this case the equations do not generate the ideal of the
modular curve, which is a curve of degree nine given by:
$s_1^2s_7+s_4^2s_1+s_7^2s_4$,
$s_1s_7^2+s_4s_1^2+s_7s_4^2-s_3^3$. The four equations above have also
four isolated points as solution: $(0:1:0:0)$ and three points
$(1:0:\om:\om^2)$ with $\om^3=1$. For
these points the equations in the $X_\al$ do not define a
curve.

There are $36$ cusps.
The relevant subgroups of $\Z/n\times\Z/n$ are the ones generated by $(1,l)$,
$0\leq l<9$,  $(0,1)$, $(3,1)$ and $(6,1)$.
The nine cusps of the last three subgroups satisfy $s_3=0$, so they are the
intersection points of the plane cubics $s_1^2s_7+s_4^2s_1+s_7^2s_4$ and
$s_1s_7^2+s_4s_1^2+s_7s_4^2$.
The other 27 cusps lie above the intersection of $s_1^2s_7+s_4^2s_1+s_7^2s_4$
with its Hessian $s_1^3+s_4^2+s_7^2-3s_1s_4s_7$.

\subex{11}
I write only the equations for the modular curve $X(11)$
using the variables $s_1$, $s_3$, $s_9$, $s_5$ and $s_4$
(in this order) \cite{\rfkf,
V.5 \S2}. The automorphism group of the curve is a group
$G_{660}$ of order $660$; I single out the transformations $U$
which is congruent to $\left({3\atop0}{0\atop4}\right)$
acting by cyclic permutation of the coordinates (namely by $s_\al
\mapsto s_{3\al{}}$), and $S$ which
multiplies $s_\al$ by $\ep^{-\ha{\al(11-\al)}}$, where $\ep=e^{2\pi
i/11}$.

Specialising the equations $F_{\al\bet\gam\del}$
yields ten equations in different eigenspaces for the induced action
of $S$. They split into two groups of five, which are permuted among
each other by the action of $U$. Therefore it suffices to write down
the following two:
$$ \macaulay
\displaylines{
s[1]s[9]3+s[9]s[5]3+s[3]s[4]3 \cr
s[1]2s[9]s[5]-s[3]s[9]2s[5]+s[1]s[3]s[4]2\;. \cr}
$$
The equations do not generate the homogeneous ideal of the curve, as
there is a zero dimensional embedded component. By subtracting the
second equation multiplied with $s_4$ from $s_1$ times the first one
obtains an expression which is divisible by $s_9$. In total one gets
five extra equations of the type:
$$\macaulay
s[3]2s[5]2+s[3]s[4]3-s[1]s[3]2s[4]+s[1]s[9]s[5]s[4]\;.
$$
Felix Klein's original way to obtain the equations starts from
the invariant
$$
\Phi=s_1^2s_3+s_3^2s_9+s_9^2s_5+s_5^2s_4+s_4^2s_1
$$
of the group $G_{660}$. Its Hessian is a quintic threefold with the
modular curve as double locus, and the fifteen $4\times4$ minors give
the equations of this curve.

\subex4
There are two equations:
$$
\eqalign{
F^0_{012}&:\quad s_1^2X_0^2-s_2^2X_1X_3+s_1^2X_2^2 \cr
F^1_{012}&:\quad s_1^2X_1^2-s_2^2X_0X_2+s_1^2X_3^2\;. \cr}
$$
The second equation is $F_{0123}$; to form the first one needs
non-integral indices. Note that $X_0(u_0)=0$,  $X_1^2(u_0)=-is_1^2$,
$X_2^2(u_0)=-s_2^2$  and $X_3^2(u_0)=is_1^2$.

\subex6
The equations $F^0_{ijk}$ are
$$
\displaylines{
s_1s_3X_0^2-s_2^2X_1X_5+s_1^2X_2X_4 \cr
s_2^2X_0^2-s_3^2X_1X_5+s_1^2X_3^2 \cr
s_1^2X_0^2-s_3^2X_2X_4+s_2^2X_3^2 \cr
s_1^2X_1X_5-s_2^2X_2X_4+s_1s_3X_3^2\;. \cr}
$$
For $2h$ odd  there is only one equation. For $h=\ha1$ one has
$$
s_1s_2X_0X_1-s_2s_3X_2X_5+s_1s_2X_3X_4\;,
$$
which can be divided by $s_2$ to give:
$$
s_1X_0X_1-s_3X_2X_5+s_1X_3X_4\;.
$$
The relation between the coefficients $s_\al$ is:
$$
s_1^4-s_2^4+s_1s_3^3\;.
$$
Remark that only even powers of $s_2$ occur (after the division).

\subex8
Even as polynomials in the $X_\al$ and $s_\al$ the
ten equations $F^0_{ijk}$ are not linearly independent.
The substitution
$(h;i,j,j,k)\mapsto (h+4;i,j,k)$ increases all indices $\al$
of the $X_\al$ by $4$; it sends e.g. $F^0_{012}$ to $F^0_{234}$.
The automorphism $(s_1,s_2,s_3,s_4)\mapsto (-s_3,-s_2,s_1,s_4)$ of order
four of the base curve can be extended to the total space by
$(X_0,X_1,X_2,X_3,X_4,X_5,X_6,X_7)\mapsto (X_0,X_5,X_2,X_7,X_4,
X_1,X_6,X_3)$. The different types of equations are:
$$
\macaulay
\catcode `\^^M\active \let ^^M\cr%
\displaylines{%
s[1]s[3]X[0]2-s[2]2X[1]X[7]+s[1]2X[2]X[6]
s[2]s[4]X[0]2-s[3]2X[1]X[7]+s[1]2X[3]X[5]
s[3]2X[0]2-s[4]2X[1]X[7]+s[1]2X[4]2
s[2]2X[0]2-s[4]2X[2]X[6]+s[2]2X[4]2
s[1]s[3]X[1]X[7]-s[2]s[4]X[2]X[6]+s[1]s[3]X[3]X[5]\;.
}$$
With $2h=1$ one obtains (out of four):
$$
\macaulay
\catcode `\^^M\active \let ^^M\cr%
\displaylines{%
s[1]s[4]X[0]X[1]-s[2]s[3]X[2]X[7]+s[1]s[2]X[3]X[6]
s[2]s[3]X[0]X[1]-s[3]s[4]X[2]X[7]+s[1]s[2]X[4]X[5]\;.
}$$
The curve defined by the relations between the coefficients
$s_\al$ is determinantal, given by the minors of the matrix:
$$
\macaulay\pmatrix{
s[3]2-s[1]2  &s[2]s[4]     &s[1]s[3] \cr
s[4]2        &s[1]2+s[3]2  &s[2]2    \cr}.
$$
It has four ordinary double points, at the coordinate vertices $e_2$
and $e_4$, and at $s_2=s_4=s_1^2+s_3^2=0$. In particular, if $s_4=1$
then $s_2=s_3^4-s_1^4$ and $s_1s_3-(s_3^4-s_1^4)^2(s_3^2-s_1^2)=0$.
The curve over the branch with tangent $s_1=0$ specialises for $s_3
\to 0$ to the $8$-gon joining the vertices in order, while the other
branch gives the one joining $e_i$ to $e_{i+3}$.

\beginsection The degeneration

\subsec
I concentrate on the case $n$ odd because the
relation to the universal curve is then more direct,
but the final formulas will also hold for even $n$.

The $n$-gon  with equations $X_\al X_\bet=0$ for
$|\al-\bet|>1$ is the fibre of the universal family  over the
cusp $P$ with coordinates $s_{n-1\over2}=1$,
$s_\al=0$ for $1\leq\al<(n-1)/2$; it corresponds to
$\tau=i\infty$.  In \cite{\rfdr, VII.4} a Tate curve with $n$
sides is constructed;  it is a generalised elliptic curve with
level $n$ structure, defined over the unit disc $D\subset \C$ with coordinate
$q^{1\over n}$, and it provides an isomorphism between the germs
$(D,0)$ and $(X(n),P)$. A modification of this construction leads to
the versal  deformation of the $n$-gon as abstract curve of genus
one, i.e., disregarding the group structure. As to infinitesimal
deformations, for an $n$-gon $Y$ one computes that $\dim T^1_Y=n$.

\roep Construction.
Let $S$ be a smooth analytic germ and consider $n$ sections
$t_0,\dots,t_{n-1}\in\Gamma(S,\sier S)$. Define $t_i$ for all $i\in\Z$ by
$t_{i+n}=t_i$. Set $T:=\prod_{i=0}^{i=n-1} t_i$.
Construct a space $\wl Y$ from charts  $(U_i)_{i\in\Z}$
with $U_i\subset S\times \C^2_i$ given by $x_iy_i-t_i=0$;
here $\C^2_i$ is a copy of $\C^2$ with coordinates $(x_i,y_i)$.
The glueing from $U_i$ to $U_{i+1}$ is determined
by $x_{i+1}y_i=1$. This makes $V_i:=U_i\cap U_{i-1}$
into $S\times \C^*$ with coordinate ring $\sier S \{x_i,x_i^{-1}\} =
\sier S \{y_{i-1},y_{i-1}^{-1}\}$, as one has $y_i=t_ix_i^{-1}$ and
$x_{i-1}=t_{i-1}y_{i-1}^{-1}$. Outside $\prod t_i=0$ all $U_i$ are identified
and one has $x_1=t_0^{-1}x_0$, $x_2=(t_0t_1)^{-1}x_0$, etc.

Let $g$ be the section of $V_n$ with $x_n(g)=1$, so $x_0(g)=T$.
Multiplication by $g$ is a well-defined operation on the union $Y$ of all
$V_i$,
and extends to $\wl Y$; it is given by  $g\colon V_i\to V_{i+n}$, $x_{i+n}(a)=
x_i(a)$ for every section of $Y\to S$.
The quotient $\wl Y / g^{\Z}$ is a family of curves of genus one.

\subsec
If all $t_i$ are equal to the same section $t$, then a group structure
can be introduced on $Y$ by  putting
$V_i\times V_j\to V_{i+j}\colon x_{i+j}(ab)=x_i(a)x_j(b)$. It extends to an
action of $Y$ on $\wl Y$.
Let $g$ be a section of $Y$, contained in a $V_n$ with $n\neq0$.
Multiplication by any section $g$  of $Y$, contained in a $V_n$ with $n\neq0$,
defines a $\Z$-action on $\wl Y$; the quotient
$\wl Y / g^{\Z}$ is a {\sl Tate curve\/}.
In particular, if $t=0$ and $g$ is the section $x_n=1$ of $V_n$, then
$\wl Y/g^{\Z}$ is the standard $n$-gon over $S$.

For $S$ the unit disc $D\subset \C$ with coordinate
$q^{1\over n}$, section $t=q^{1\over n}$ and section $g$ of $V_n$
with $x_n(g)=1$ (so $x_0(g)=q$), denoted by $q$ again,
the quotient $\wl Y/q^{\Z}$
is the Tate curve with $n$ sides. By the modular interpretation the $s_\al$ can
be
written as functions of the variable $q^{1\over n}$. An
explicit $q$-expansion can be derived from formula (10.6) in
\cite{\rfvel}, or from formula (3) p.~281 in \cite{\rfkf}. Here
I only note that $s_a$ is  in first approximation proportional
to $(-1)^{-\al{}} (q^{1\over n})^{-\ha{\al(n-\al)}}$.

\subproclaim Theorem.
\xdef\labelth{\the\secno.\the\subsecno}
The versal deformation of the $n$-gon has a smooth base space
with deformation parameters $t_i$; it is given by
the equations
$$
\displaylines{\qquad
F^{h}_{ijk}=
s_{k+j}s_{k-j}
\left(\prod_{m=i}^{j-1}T_{h-m}^{h+m}\right)X_{h+i}X_{h-i}-
s_{k+i}s_{k-i}X_{h+j}X_{h-j}
\hfill\cr\hfill{}+
s_{j+i}s_{j-i}
\left(\prod_{m=j+1}^{k}T_{h+m}^{h-m}
\right)X_{h+k}X_{h-k}\;.\quad\cr}
$$
Here $T_i^j :=  \prod _{k=i}^j t_k$, and $T:=  \prod _{k=0}^{n-1} t_k$.
Furthermore \/$h \in {\Bbb Z}/n{\Bbb Z}$\/ for odd $n$, and
\/$h \in {1\over 2}{\Bbb Z}/{1\over2}n{\Bbb Z}$ for even $n$, and\/
$0\leq i<j<k \leq n/2$; for even $n$ the numbers $i$, $j$ and $k$ are
integral if and only if $h$ is integral.
The equations $F^h_{j-1,j,j+1}$ form a A minimal set of generators
of the ideal.
The coefficients $s_\al$ in the
equations are holomorphic functions of $T$ with
$s_\al(0)=(-1)^{n(\al+1)}$ and  they satisfy:
$$
-s_{k+i}s_{k-i}s_{l+j}s_{l-j}+
s_{l+i}s_{l-i}s_{k+j}s_{k-j}+
s_{j+i}s_{j-i}s_{l+k}s_{l-k}T^{k-j}=0\;,
$$
for\/ $0\leq i<j<k<l \leq n/2$.

\roep Proof.
I will not express the functions $X_\al$ in terms of the coordinates
on  the space $Y$, so the proof  will be indirect: I
shall establish that the given equations define  a deformation
of the $n$-gon by studying the syzygies.

As explanation of the equations I offer the circumstance that for
fixed $T\neq0$ they are really the same as Klein and Fricke's
equations: I  describe an explicit coordinate transformation. For
odd $n$ it involves the $n$th root of each $t_i$, whereas for even $n$
one needs the $2n$th root. I  write $t_i=\tau_i^{n}$ for all $n$,
so for even $n$ there is the additional indeterminacy of a
square root. Consider the transformation which replaces $X_i$ by
$\prod_j\tau^{a_{ij}}X_i$, and $s_i$ by $(\prod_j\tau_j)^{b_i}s_i$ with
$a_{ij}=(j-i)(n-j+i)/2$ for $i\leq j \leq i+n$;
the exponents $a_{ij}$  are solutions to the equations:
$$
\eqalign{2a_{ij}&=a_{i,j-1}+a_{i,j+1}+1\;, \qquad j\neq i\;,\cr
n&=a_{i,i-1}+a_{i,i+1}+1\;, \cr
a_{i,i-j}&=a_{i,i+j}\;, \cr
a_{ii}&=0\;. \cr}
$$
The $b_i$ satisfy similar equations: $b_i=(i-1)(n-i-1)/2$ for $0< i<n$.
Note that the exponent of $q^{1\over n}$ in the $q$-expansion of $s_i$ given
above is $-i(n-i)/2$, which differs from $-b_i$ by $(n-1)/2$,
independently of $i$.

This coordinate change transforms the equations to multiples of the
elliptic curve equations. I check here only the equations for the
parameters. The exponent of $\prod \tau_i$ as coefficient
of the monomial $s_{k+i}s_{k-i}s_{l+j}s_{l-j}$ is:
$$
\displaylines{
\qquad
\textstyle {(k+i-1)(n-i-k-1)\over2}+{(k-i-1)(n-k+i-1)\over2}+
{(l+i-1)(n-i-l-1)\over2}+{(l-i-1)(n-l+i-1)\over2}
\hfill\cr\hfill{}=
(k+l-2)n-k^2-i^2-l^2-j^2+2\;, \qquad\cr}
$$
which  by symmetry is  the same as for
$s_{l+i}s_{l-i}s_{k+j}s_{k-j}$; due to the extra term $T^{k-j}$
the third summand gives:
$$
((j+l-2)n-j^2-i^2-l^2-k^2+2)+n(k-j)=(k+l-2)n-k^2-i^2-l^2-j^2+2\;.
$$
Note that for even $n$ this is an integral power of $\prod \tau_i$.

The $q$-expansion mentioned above shows that solutions to the
$s_\al$ equations  exist, which are  power series in $T$.

\proclaim Lemma 1.
Let $0\leq i<j<k<l\leq n/2$. The following relations hold:
$$
\def\dp{\displaystyle}
\vcenter{\openup 1\jot
\halign{$\dp#$&$\dp#$&${}#$\cr
R^{h;i}_{jkl}\colon\quad&
s_{l+i}s_{l-i}F^{h}_{ijk}-s_{k+i}s_{k-i}F^{h}_{ijl}
+s_{j+i}s_{j-i}
\left(\prod_{m=j+1}^{k}T_{h+m}^{h-m}\right)F^{h}_{ikl}&=0\cr
R^{h;j}_{ikl}\colon\quad&
s_{l+j}s_{l-j}F^{h}_{ijk}-s_{k+j}s_{k-j}F^{h}_{ijl}
+s_{j+i}s_{j-i}
\left(\prod_{m=j+1}^{k}T_{h+m}^{h-m}\right)F^{h}_{jkl}&=0\cr
R^{h;k}_{ijl}\colon\quad&
s_{l+k}s_{l-k}\left(\prod_{m=j}^{k-1}T_{h-m}^{h+m}\right)F^{h}_{ijk}
-s_{k+j}s_{k-j}F^{h}_{ikl}
+s_{k+i}s_{k-i}F^{h}_{jkl}&=0\cr
R^{h;l}_{ijk}\colon\quad&
s_{l+k}s_{l-k}\left(\prod_{m=j}^{k-1}T_{h-m}^{h+m}\right)F^{h}_{ijl}
-s_{l+j}s_{l-j}F^{h}_{ikl}
+s_{l+i}s_{l-i}F^{h}_{jkl}&=0\;.\cr}}
$$

\noindent The proof consists of a simple verification. The
assumption $0\leq i<j<k<l\leq n/2$ and the division into four
cases is needed to  be able to specify the terms  involving the
$t_\al$.


The four-term relations are best formulated in terms
of equations of the form
$F_{\al_0\al_1\al_2\al_3}$,  The precise form of the equation
$F_{\al_0\al_1\al_2\al_3}$ depends on the order of the six
points $\al_0+\al_1$, $\al_2+\al_3$, $\al_0+\al_2$, $\al_3+\al_1$,
$\al_0+\al_3$, $\al_1+\al_2$. If $\al_0<\al_1<\al_2<\al_3$, then
$X_{\al_0+\al_2}X_{\al_3+\al_1}$ is the `middle monomial', the one
occurring without $t_i$-term, as $\al_0+\al_1$ and $\al_2+\al_3$ lie
between $\al_1+\al_3$ and $\al_0+\al_2$, whereas $\al_0+\al_3$ and
$\al_1+\al_2$ lie on the other side.  As I really work modulo
$n$\/ I can always arrange for a given index to be the smallest.

\proclaim Lemma 2.
Suppose $0\leq\al_0<\al_1<\al_2<\al_3<\al_4$.
Let $D_2$ be the product of factors $T_i^j$, occurring in
$F_{\al_1\al_2\al_0\al_4}$ as coefficient of
$X_{\al_1+\al_0}X_{\al_2+\al_4}$, and $C_2$ that of
$X_{\al_0+\al_4}X_{\al_1+\al_3}$ in
$F_{\al_1\al_0\al_3\al_4}$.
Then
$$
D_2s_{\al_0-\al_1}X_{\al_0+\al_1}F_{\al_0\al_2\al_3\al_4} +
s_{\al_0-\al_2}X_{\al_0+\al_2}F_{\al_1\al_0\al_3\al_4}
+s_{\al_0-\al_3}X_{\al_0+\al_3}F_{\al_1\al_2\al_0\al_4} +
C_2 s_{\al_0-\al_4}X_{\al_0+\al_4}F_{\al_1\al_2\al_3\al_0}
=0\;.
$$

\roep Proof.
\def\al{\alpha}
Write $B_1$, \dots, $E_2$ for the product of the factors $T_i^j$
occurring in the equations involved.
The same substitution as before shows the existence of a relation of
the form
$$
\displaylines{\quad
A_1s_{\al-\bet}
X_{\al+\bet}F_{\al\gam\del\ep} +
A_2s_{\al-\gam}X_{\al+\gam}F_{\bet\al\del\ep}
+A_3s_{\al-\del}X_{\al+\del}F_{\bet\gam\al\ep}+
A_4s_{\al-\ep}X_{\al+\ep}F_{\bet\gam\del\al}=
\hfill \cr%
\quad A_1s_{\al-\bet}
X_{\al+\bet}
\cdot\Bigl(s_{\al-\gam}s_{\del-\ep}
B_1X_{\al+\gam}X_{\del+\ep} +
s_{\al-\del}s_{\ep-\gam}X_{\al+\del}X_{\ep+\gam} +
B_2s_{\al-\ep}s_{\gam-\del}
X_{\al+\ep}X_{\gam+\del}\Bigr) \hfill\cr
\quad {}+
A_2s_{\al-\gam}X_{\al+\gam}
\cdot\Bigl(C_1s_{\bet-\al}s_{\del-\ep}X_{\bet+\al}X_{\del+\ep}
+ C_2s_{\bet-\del}s_{\ep-\al}
X_{\bet+\del}X_{\ep+\al} +
s_{\bet-\ep}s_{\al-\del}X_{\bet+\ep}X_{\al+\del}\Bigr)
\hfill\cr
\quad {}+
A_3s_{\al-\del}X_{\al+\del}
\cdot\Bigl(D_1s_{\bet-\gam}s_{\al-\ep}
X_{\bet+\gam}X_{\al+\ep}+
D_2s_{\bet-\al}s_{\ep-\gam}
X_{\bet+\al}X_{\ep+\gam} +
s_{\bet-\ep}s_{\gam-\al}X_{\bet+\ep}X_{\gam+\al}\Bigr)
\hfill \cr
\quad {}+
A_4s_{\al-\ep}
X_{\al+\ep}
\cdot\Bigl(E_1s_{\bet-\gam}s_{\del-\al}
X_{\bet+\gam}X_{\del+\al} +
s_{\bet-\del}s_{\al-\gam}X_{\bet+\del}X_{\al+\gam} +
E_2s_{\bet-\al}s_{\gam-\del}
X_{\bet+\al}X_{\gam+\del}\Bigr),
\hfill \cr}
$$
with coefficients $A_i$ satisfying
$A_1B_1=A_2C_1$, $A_1=A_3D_2$, $A_1B_2=A_4E_2$, $A_2C_2=A_4$,
$A_2=A_3$, and $A_3D_1=A_4E_1$. These equations show that I can
take $A_2=A_3=1$, $A_1=D_2$ and $A_4=C_2$.

\proclaim Lemma 3.
The  equations $F^h_{ijk}$ generate the
ideal of the total space of the degeneration, while the  relations of
Lemma 1 and Lemma 2
generate the syzygies between them.

\roep Proof.
The relations $R^{h;j}_{ikl}$ and $R^{h;k}_{ijl}$
show that all equations $F^{h}_{ijk}$ can be expressed terms
of the equations  $F^{h}_{j-1,j,j+1}$. For odd $n$ and fixed
$h\in {\Bbb Z}/n{\Bbb Z}$ there are $(n-3)/2$ such equations,
whereas for even $n$ one has $(n-2)/2$ equations with $2h$ even and
$(n-4)/2$ with $2h$ odd. If $t_i=0$ for all $i$, these equations
generate the ideal of the $n$-gon.

The equations $F^{h}_{j-1,j,j+1}$ are linearly independent;
for fixed $h$ one can use
induction on $j$: the monomial $X_{h+j-1,h-j+1}$ occurs first in
the equation $F^{h}_{j-1,j,j+1}$ with nonzero coefficient, as
$T\neq0$; note that if $T=0$, then from some $j_0$ onwards the
smallest monomial occurring in $F^{h}_{j-1,j,j+1}$ is
$X_{h+j}X_{h-j}$.

As to the relations is suffices to lift the generators of the syzygy module
for the $n$-gon, which are of the form $(X_iX_j)X_k-(X_iX_k)X_j=0$.
I write the formula for $(i,j,k)=(2\gam, \gam+\bet,\gam+\al)$. To be
specific I suppose that $\al<\bet<\gam$.
The two equations to consider are
$$
\displaylines{\quad
F_{\al,\gam,\gam-1,\gam+1}=\hfill\cr\hfill
s_{\al-\gam+1}s_{1}T_{2\gam+1}^{\al+\gam-1}X_{\al+\gam-1}X_{2\gam+1}
-s_{\al-\gam}s_{2}X_{\al+\gam}X_{2\gam} +
s_{\al-\gam-1}s_{1}T^{2\gam-1}_{\al+\gam+1}X_{\al+\gam+1}X_{2\gam-1}
\quad\cr}
$$
and
$$
\displaylines{\quad
F_{\gam,\bet,\gam-1,\gam+1}=\hfill\cr\hfill
s_{\gam-\bet+1}s_{1}T_{2\gam+1}^{\bet+\gam-1}X_{\bet+\gam-1}X_{2\gam+1} -
s_{\gam-\bet}s_{2}X_{\bet+\gam}X_{2\gam} +
s_{\gam-\bet-1}s_{1}T^{2\gam-1}_{\bet+\gam+1}X_{\bet+\gam+1}X_{2\gam-1}\;.
\quad\cr}
$$
They occur in the relation
$$
\displaylines{\qquad
s_{\gam-\al}X_{\gam+\al}F_{\gam,\bet,\gam-1,\gam+1} +
s_{\gam-\bet}X_{\gam+\bet}F_{\al,\gam,\gam-1,\gam+1}
\hfill\cr\hfill{}+
s_{1}X_{2\gam-1}T^{2\gam-1}_{\bet+\gam+1}F_{\al,\bet,\gam,\gam+1} -
s_{1}X_{2\gam+1}T_{2\gam+1}^{\al+\gam-1}F_{\al,\bet,\gam-1,\gam}
=0\;.\qquad\cr}
$$

\proclaim Lemma 4.
The deformation  in the $t_\al$-direction smoothes the double
point of the $n$-gon at the vertex $e_\al$ of the coordinate
simplex, whereas the curve has a node at $e_\al$ if $t_\al=0$.

\roep Proof.
Consider the standard affine charts in $\P^{n-1}$. It suffices to look at
$X_0=1$. Let $x_\al=X_\al/X_0$. All coordinates $x_\al$ except $x_1$ and
$x_{n-1}$ can be eliminated successively. I start with the equation
$$
F^{-1}_{012}:\qquad s_2^2x_{n-2}=
s_3s_1t_{n-1}x_{n-1}^2+s_1^2T_1^{n-3}x_1x_{n-3}\;.
$$
Using it the equation $F^{-\ha3}_{\ha1\ha3\ha5}$ gives
$$
(s_{3}s_{2}-s_{4}s_1^3s_2^{-2}T_1^{n-1}x_1x_{n-1})x_{n-3}=
s_4s^3s_1^2s_2^{-2}t_{n-2}t_{n-1}^2x_{n-1}^3
+s_{2}s_{1}T_1^{n-4}x_1x_{n-4}\;.
$$
In the neighbourhood of the origin
the coefficient of $x_{n-3}$ is a unit so the equation can be used
to eliminate $x_{n-3}$.
By induction
$$
x_{n-i}= ({\rm unit})\cdot t_{n-i+1}t_{n-i+2}^2\cdots
t_{n-1}^{i-1}x_{n-1}^i+ ({\rm unit})\cdot T_1^{n-i-1}x_1x_{n-i-1\;}.
$$
In particular,
$$
x_2= ({\rm unit})\cdot t_{3}t_{4}^2\cdots
t_{n-1}^{n-3}x_{n-1}^{n-2}+ ({\rm unit})\cdot t_1x_1^2\;.
$$
Going back one finds
$$
x_i= ({\rm unit})\cdot t_{i+1}t_{i+2}^2\cdots
t_{n-1}^{n-i-1}x_{n-1}^{n-i}+
({\rm unit})\cdot t_{i-1}t_{i-2}^2\cdots t_1^{i-1}x_1^i\;.
$$
Finally the equation $F^0_{012}$ gives
$$
({\rm unit})\cdot x_1x_{n-1}- ({\rm unit})\cdot t_{n-2}t_{n-3}^2\cdots
t_1^{n-2}x_1^n-({\rm unit})\cdot t_{2}t_{3}^2\cdots
t_{n-1}^{n-2}x_{n-1}^n
= t_0\;,
$$
which proves the claims.
\qed

\subex6
I write two of the equations $F^0_{ijk}$ (the other two are
obtained by replacing $X_\al$, $t_\al$ by $X_{\al+3}$ and
$t_{\al+3}$):
$$
\displaylines{
s_1s_3t_0X_0^2-s_2^2X_1X_5+s_1^2t_2t_3t_4X_2X_4 \cr
s_2^2t_0X_0^2-s_3^2X_1X_5+s_1^2t_2t_3^2t_4X_3^2\;. \cr
}
$$
Furthermore  $F^{\ha1}_{{1\over2}{3\over2}{5\over2}}=
s_2(s_1t_0t_1X_0X_1-s_3X_2X_5+s_1t_3t_4X_3X_4)$.
The relation between the coefficients is $s_1^4T-s_2^4+s_1s_3^3$.
A solution is $(s_1,s_2^2,s_3)=(1,\sqrt{1+T},1)$.
Ruud Pellikaan obtained this formula (in the cusp case, with
$t_i=X^{a_i-2}$)  by direct computations with syzygies.

\subex9
The  equations on the coefficients become:
$$
\displaylines{
s_1^3s_4T-s_1s_3^3-s_7^3s_4  \cr
s_1^3s_7T^2+s_3^3s_4-s_4^3s_7  \cr
s_7^3s_1T-s_3^3s_7-s_4^3s_1  \cr
s_1^2s_3s_7T+s_7^2s_3s_4+s_4^2s_1s_3\;.  \cr}
$$
As the $s_\al$ are nonzero for small $T$ the equations reduce to:
$s_1^2s_7T+s_4^2s_1+s_7^2s_4$ and
$s_1s_7^2T+s_4s_1^2T+s_7s_4^2-s_3^3$.
As solution  in closed form  one can
take   $s_1=s_7=1$,  $s_4=-\ha1(1+\sqrt{1-4T})$, and
$s_3=\left(\ha1(1-T)(1+\sqrt{1-4T})\right)^{1/3}$.

\beginsection Cusp singularities

\subsec
A cusp singularity is a normal isolated two-dimensional singularity,
whose minimal resolution has a cycle of rational
curves as exceptional divisor. The analytical type of the singularity
is completely determined by the self-intersections $-b_i$ of these
curves. If the exceptional divisor is irreducible, it is a nodal
curve.

The name cusp comes from the connection with
Hilbert modular forms;
the standard reference is \cite{\rfhir}; see also
\cite{\rfbeht, \rfgeer}.
Let $K=\Bbb Q(\sqrt D)$ be a real quadratic field.
Two lattices $M_1$, $M_2$ in $K$ are strictly equivalent if $M_2=\lambda M_1$
for some totally positive $\lambda\in K$ (i.e., both $\lambda$ and
its conjugate $\lambda'$ are positive).  A lattice $M$ is always
strictly equivalent to one of the form ${\Bbb Z}\oplus{\Bbb Z}w$
with $0<w'<1<w$. Let $U^+_M$ be the group of totally positive
elements $\ep$ of $K$, satisfying $\ep M=M$, and let $V$ be a
subgroup of finite index $s$ in $U^+_M$. The semi-direct product
$G(M,V)=M\rtimes V$ acts  freely and properly discontinuously on
the product of upper half spaces
${\frak H}^2$. The quotient $G(M,V)\backslash {\frak H}^2$ can be
compactified by one cusp $\infty$ to a normal complex space
$X(M,V)$. The number $w$ has a periodic
continued fraction expansion $w=[[b_1,\dots,b_r]]$ with
primitive period of length $r/s$.
The numbers $b_i$ are the same as before, minus the self-intersection of
the exceptional curves on the minimal resolution $\wtx$, except in the case
$r=1$; then  $-b_1$ is the degree of the normal bundle of the composed map
$\Bbb P^1\to E_1\subset \wtx$.

The resolution can be constructed by toroidal methods. To this end
one considers the embedding
of $M_+$, the totally positive elements
of $M$, in ${\Bbb R}^2_+$, given by the two embeddings of $K$ into $\Bbb R$.
The lattice points on the boundary of the convex hull of $M_+$ define
a fan and with it one constructs the manifold $\overline Y$ containing
an infinite chain of rational curves, which comes with an embedding
of $Y:=M\backslash {\frak H}^2$, and an action of $V$. The orbit space
$\overline Y/V$ is the minimal resolution of $X(M,V)$. This
construction is similar to that of the Tate curve in the
previous section.

The dual lattice $M^*=\{\mu \in K\mid \Tr(\al\mu)\in {\Bbb Z}
\hbox{\rm~for~all~}\al\in M\}$ is strictly equivalent to
$M={\Bbb Z}\oplus{\Bbb Z}w^*$ with $w_*=(2-w')/(1-w')$; here
$\Tr(\al)=\al+\al'$. The
continued fraction $w_*=[[a_1,\dots,a_n]]$ (again with the primitive
period repeated $s$ times) is the dual of
$[[b_1,\dots,b_r]]$ in reverse order: the sequence $(b_1,\dots,b_r)$
starts with $a_n-3$ twos.

The local ring of the cusp consists of certain Fourier series, of the  form
$$
f(z)=\sum_{\mu\in M^*\cup\{0\}} a_\mu \exp 2\pi i(\mu z_1 +\mu' z_2)\;,
$$
with $a_\mu=a_{\ep\mu}$ for all $\ep \in V$. Abbreviate
$e(\mu z)= \exp 2\pi i(\mu z_1 +\mu' z_2)$.
Consider the embedding of $M^*$ in ${\Bbb R}^2$  defined by
$i(\mu)=(\mu,\mu')$. Let $i(A_i)$ be the lattice points on the boundary of the
convex hull of $i(M^*_+)$. Put $A_0=1$,  $A_{-1}=w_*$ and define
$w_k:=A_{k-1}/A_k$. One has $A_{k-1}+A_{k+1}=a_kA_k$, or
equivalently $w_k=a_k-1/w_{k+1}$, so $w_k=
[a_k,a_{k+1},\dots]$. Then $A_n$ is the unique generator $\ep_1$
of $V$ with  $0<\ep_1<1$. For $\mu \in M^*$ define `Poincar\'e
series'  $F_\mu(z):=\sum_{\ep\in V}e(\ep\mu z)$. If $n\geq3$, the
$n$ functions $F_{A_1}$, \dots, $F_{A_n}$ generate the local ring
of the cusp, and define an embedding in $({\Bbb C}^n,0)$; for
$n=1,2$ one has to take some extra generators (see
\cite{\rfcoh}).  The exact relations between these generators are
not known. Even in the  hypersurface case ($n\leq 3$) the
equation is not polynomial, but its Newton diagram has the
standard form (for an explicit example see \cite{\rfcoh}).

\subproclaim Proposition.
The specialisation $t_i=X_i^{a_i-2}$ in the equations for the deformation of
the  $n$-gon in Theorem \/{\rm (\labelth)} gives rise to defining equations of
a cusp singularity with $w_*=[[a_1,\dots,a_n]]$.

\roep Proof.
The equations of (\labelth) describe a deformation of the cone over the
$n$-gon over a smooth base space. The specialisation equations
$t_i-X_i^{a_i-2}$ form a regular sequence, so the  resulting ideal is the
ideal of a surface singularity. Equations of the tangent cone are obtained by
putting $t_i=0$ if  $a_i>2$ (and $t_i=1$ if $a_i-2=0$). As at least one
$a_i>2$, the product $T=\prod t_i$ vanishes. The projectivised tangent cone
consists of a cycle of rational curves, with as many double points as the
number of $i$ with $a_i>2$.

Singular points on the first blow-up can only occur in double points of the
exceptional locus. Suppose $a_0>2$. Let $X_0=x_0$, $X_i=x_0x_i$ be the
description of the blow-up in one affine chart. The strict transform of the
surface is given by the same equations as in the proof of Lemma 4, with
$t_i=x_0^{a_i-2}x_i^{a_i-2}$ for $i\neq0$ and $t_0=x_0^{a_0-2}$.
It is smooth at the origin if $a_0=3$, and otherwise there is a singularity
of type $A_{a_0-3}$. As cusp singularities are characterised by their
exceptional divisor, this computation shows that the surface singularity is
indeed a cusp associated with $w_*=[[a_1,\dots,a_n]]$.
\qed
\endroep

By construction the cusp comes with a deformation to a
simple elliptic singularity of multiplicity $n$.
Such a singularity is in fact the general fibre over a
$\sum(a_i-2)$-dimensional base space:
introduce deformation parameters $v_i^{(j)}$, $j=1,\dots,a_i-2$, and
put
$$
t_i=X_i^{a_i-2}+v_i^{(1)}X_i^{a_i-3}+\dots+v_i^{(a_i-3)}X_i+v_i^{(a_i-2)}\;.
$$
This family has codimension $n$ in $T^1$, as $\dim T^1=\sum(a_i-1)$
\cite{\rfbeht}. For $n>9$ it forms an irreducible component of the
deformation space: simple elliptic singularities of multiplicity $n>9$ have
only equisingular unobstructed deformations, of codimension $n$.

This deformation is similar to the Artin
component of cyclic quotient singularities. There
the additional deformations are obtained by adding an extra
summand $v_i^{(a_i-1)}/X_i$ to each $t_i$
\cite\rfarn. In the case at hand it
is difficult to make sense of such a formula, as the $s_i$ depend
on the $t_i$; for $n=6$ this idea works, as the formulas in the next section
show.

\subsec
Degenerate cusp  are a class of weakly normal, non-isolated surface
singularities, first introduced and described by Shepherd-Barron
\cite{\rfsb}.  For such singularities the role of the resolution is taken
over by  an improvement. By definition an {\sl improvement\/} $\pi\colon Y\to
X$ is a point modification, such that $Y$ is weakly normal, the singular locus
$S$ of $Y$ is the strict transform  of $\Sing (X)$, while $S$,  the
normalisation $\widehat Y$ of $Y$ and the   inverse image $\widehat S$ of $S$
on $\widehat Y$ are all smooth \cite{\rfjs}. Define $Y$ to be {\sl
weakly minimal\/}, if $\widehat Y$  contains only
$(-1)$-curves that intersect $\widehat S$.

\roep Definition.
A degenerate cusp singularity  is a weakly normal surface singularity,
such that a weakly minimal improvement has a cycle of rational
curves as exceptional divisor.
\endroep

The normalisation $\widehat X$ of $X$ consists of a
disjoint union of cyclic quotient singularities $X_i$ (including the case of
$A_0$, a smooth germ), and the inverse image of $\Sing(X)$ consists of
two transversal lines. To construct the improvement one takes the minimal
embedded resolution $Y_i\to X_i$ of these lines, and glues the $Y_i$ along
the strict transform of the lines to form $Y$.

One can still associate a cycle of numbers $[[a_1,\dots, a_n]]$ to
$X$. In the non-degenerate case the canonical model contains an
$A_{i-3}$-singularity
for each $i$ with $a_i\geq3$; on the `canonical improvement' of a
degenerate cusp
(which for $n\geq3$ is just the first blow-up)  an
$A_\infty$-singularity then corresponds to $a_i=\infty$.

\subsec
The same procedure as for non-degenerate cusps yields equations: put
$t_i=X^{a_i-2}$. The equations simplify because $t_i=0$, if $a_i=
\infty$, and therefore $T=\prod t_i=0$. This allows to take
$s_\al=1$ (if $\al$ is odd).  The resulting equations were
originally obtained by Ruud Pellikaan  from the equations for
cyclic quotients.

\proclaim Proposition {\rm(Pellikaan)}.
Let $X(a_1,\dots,a_n)$ be a degenerate cusp with $n>3$. The $n(n-3)/2$
equations
$$
F_{ij}:=F^{\ha{i+j}}_{\ha{j-i}-1,\ha{j-i},\ha{j-i}+1}\colon\quad
X_iX_j=X_{i+1}\left(\prod_{k=i+1}^{j-1}X^{a_k-2}\right)X_{j-1} +
       X_{j+1}\left(\prod_{k=j+1}^{i-1}X^{a_k-2}\right)X_{i-1}
$$
with $i-j\neq -1,0,1$ generate the ideal of the cusp.

For degenerate cusps the equations can be `blown down' if the
cycle $[[a_1,\dots, a_n]]$ contains a $1$. This means that I
formally write $t_i=1/X_i$, if $a_i=1$.
Suppose $a_n=\infty$. If $i\neq1,n-1$, then the equation
$F_{i-1,i+1}\colon X_{i-1}X_{i+1}-X_i$ allows to eliminate $X_i$,
and substituting $1/X_{i-1}X_{i+1}$ for $t_i$ shows that the
new equations correspond to
$[[a_1,\dots,a_{i-1}-1,a_{i+1}-1,\dots,a_{n-1}, \infty]]$. The case $a_{n-1}=1$
gives the  equation
$X_{n-1}=X_nX_{n-2}-X_1T_1^{n-3}X_{n-3}$. As $t_n=0$ the only
equations in which $t_{n-1}$ appears are the  $F_{kn}$, where the
combination $t_{n-1}X_{n-1}(=1)$ occurs. The $(n-1)(n-4)/2$
equations not involving the index $n-1$ define the cusp
$[[a_1,\dots,a_{n-2}-1,\infty]]$.

In case $n=4$ the result is a cubic equation. Blowing down now also
works for isolated cusps. Consider the equations for
$[[a_0,a_1,a_2,1]]$:
$$
\eqalign{
X_1X_3&=X_0^{a_0}+X_2^{a_2}\cr
X_0X_2&=X_1^{a_1}+X_3\;.      \cr
}
$$
Elimination of $X_3$ gives the equation of $T_{a_0,a_1+1,a_2}$,
the cusp with cycle $[[a_0-1,a_1,a_2-1]]$:
$$
X_0X_1X_2-X_1^{a_1+1}+X_0^{a_0}+X_2^{a_2}\;.
$$

In particular, equations for the cycle
$[[2,\dots,2,1,\infty]]$ blow down to  $T_{23\infty}$, the cusp
with cycle $[[\infty]]$.

\proclaim  Proposition.
Every degenerate cusp deforms to $T_{23\infty}$ and is
therefore smoothable.

\roep Proof.
I can suppose that $a_0=\infty$. As noted above $t_{n-1}$ occurs
only in the combination $t_{n-1}X_{n-1}$ in the equations $F_{ij}$;
but note that this does not hold for general $F^h_{ijk}$.
The relations between the $F_{ij}$ can easily be computed directly
to be:
$$
\displaylines{

\matrix{
X_kF_{ij}-X_jF_{ik}+X_{i+1}T_{i+1}^{j-1}F_{j-1,k}=0\cr
X_iF_{jk}-X_jF_{ik}+X_{k-1}T_{j+1}^{k-1}F_{i,j+1}=0\cr}
\qquad {\rm if}\quad 0<i<j<k\leq n-1\cr
\noalign{\vskip 3pt}
X_iF_{0j}-X_jF_{0i}-X_1T_1^{i-1}F_{i-1,j}+X_{n-1}T_{j+1}^{n-1}F_{i,j+1}=0\cr
X_iF_{0j}-X_0F_{ij}+X_{n-1}T_{j+1}^{n-1}F_{i,j+1}
+X_{n-1}T_{i+1}^{n-1}F_{i+1,j}=0\;.
\cr}
$$
Again $t_{n-1}$ only occurs together with $X_{n-1}$. The required
deformation can now be defined by putting $t_0=0$,
$t_i=X_i^{a_i-2}+\ep$ for $1\leq i\leq n-2$ and
$t_{n-1}=X_{n-1}^{a_{n-1}-2}+\ep /X_{n-1}$.

\beginsection Smoothings of simple elliptic singularities

In this section I use the symmetric equations to describe
deformations of simple elliptic singularities. The dimension of
$T^1$, the vector space of infinitesimal deformations, is $n+1$
for a singularity of multiplicity $n$; the graded parts  satisfy
$\dim T^1(-1) =n $ and $\dim T^1(0) =1$. I have
not been able  to find the general formula for the perturbation
of the equations $F^h_{ijk}$. Therefore I restrict myself to
$n\leq 9$, the cases in which cones over elliptic curves of
degree $n$ are smoothable. I only give the results of my computations.

\subex6
For this multiplicity the versal deformation of cusps can be given.
It is known that the base is up
to a smooth factor  independent of the
specific cusp; this circumstance enables me to give
a general formula.

I introduce deformation parameters $v_i(:=v_i^{(a_i-1)})$, while
the others are implicit in the formula
$$
t_i=X_i^{a_i-2}+v_i^{(1)}X_i^{a_i-3}+\dots+v_i^{(a_i-3)}X_i+v_i^{(a_i-2)}.
$$
These deformation parameters also enter the equations through
the formula $s_1^4T-s_2^4+s_1s_3^3$, contrary to  the $v_i$.
The deformation of the simple elliptic singularity is obtained
by putting $t_i=1$.

I describe the equations $F^0_{123}$,
$F^0_{013}$ and  $F^{1/2}_{1/2,3/2,5/2}=F_{0124}$:
$$\macaulay
\displaylines{ \quad
s[1]2t[5]t[0]t[1]X[1]X[5]-s[2]2X[2]X[4]+s[1]s[3]t[3]X[3]2
 \hfill\cr \qquad  {}
    +s[3]s[2]2v[3]X[3]+s[1]s[2]2t[0](v[5]t[1]X[1]+v[1]t[5]X[5])
         +s[1]2s[3]t[1]t[0]t[5]v[0]X[0]
        -s[1]3t[0]t[1]t[3]t[5](v[2]t[4]X[4]+v[4]t[2]X[2])
                 \hfill\cr \qquad  {}
      +s[2]4t[0]v[1]v[5]+s[1]s[3]s[2]2t[1]t[5]v[0]2
         -s[1]2s[2]2t[0]t[1]t[3]t[5]v[2]v[4]
\hfill
\cr
}
$$
$$\macaulay
\displaylines{
\quad
s[2]2t[0]X[0]2-s[3]2X[1]X[5]+s[1]2t[2]t[3]2t[4]X[3]2
 \hfill\cr \qquad  {}
         +(s[3]3+2s[1]3T)v[0]X[0]
         +2s[1]s[2]2X[3]t[2]t[3]t[4]v[3]
                 +s[1]s[3]2t[3](v[2]t[4]X[4]+v[4]t[2]X[2])
         \hfill\cr \qquad  {}
         +s[3]2s[2]2t[3]v[2]v[4]+s[2]4t[2]t[4]v[3]2
         +s[1]2s[2]2t[1]t[2]t[3]t[4]t[5]v[0]2
\hfill
\cr
}
$$
$$\macaulay
\displaylines{
\quad
s[1]t[0]t[1]X[0]X[1]-s[3]X[2]X[5]+s[1]t[3]t[4]X[3]X[4]
 \hfill\cr \qquad {}
         +s[2]2(v[3]t[4]X[4]+v[1]t[0]X[0]+v[0]t[1]X[1]+v[4]t[3]X[3])
                 +s[1]2s[3]t[0]t[1]t[3]t[4]v[2]v[5]
         +(1/2s[3]3+s[1]3T)(v[0]v[1]+v[3]v[4])
\hfill
\cr
}
$$
The base space is (up to a smooth factor) the  cone over
$\P^1\times\P^2$:
$$\macaulay
\pmatrix{v[0] & v[2] & v[4] \cr
         v[3] & v[5] & v[1] \cr}
$$
The infinitesimal deformation in the $v_0$ direction can be
obtained by formally  putting $t_0=X_0^{a_0-2}+v_0s_2^2/s_1X_0$,
using the other equations and making coordinate transformations
involving the `unit' $1+s_1^3t_1t_2t_3t_4t_5v_0/4s_2^2X_0$.

\subex7
The equations of the total space are not very
enlightening, but I give them for completeness.
The equation $F^0_{012}$ is changed into:
$$
\macaulay
\displaylines{
\quad s[1]s[4]X[0]2+s[2]2X[1]X[6]-s[1]2X[2]X[5]
+(s[2]2s[4]-s[1]3)s[1]s[2]s[4]v[0]X[0] \hfill\cr
\quad  {}+2(s[1]2s[2]-s[4]3)s[1]s[2]2(v[6]X[1]+v[1]X[6])
  -(s[4]3-s[1]2s[2])s[1]3(v[5]X[2]+v[2]X[5])
  -3s[1]3s[2]2s[4](v[4]X[3]+v[3]X[4]) \hfill\cr
\quad  {}+2(s[2]s[4]5+s[1]5s[4]-2s[1]s[2]5)s[1]2s[2]s[4]v[1]v[6]
  +(s[4]5+3s[1]2s[2]s[4]2-7s[1]s[2]4)s[4]s[1]4v[2]v[5] \hfill\cr
\quad {}
  +(4s[1]s[2]4-3s[1]2s[2]s[4]2+2s[4]5)s[1]2s[2]s[4]2v[3]v[4] \hfill\cr}
$$
and $F^0_{123}$ into:
$$\macaulay
\displaylines{
\quad s[1]s[2]X[1]X[6]+s[2]s[4]X[2]X[5]+s[1]s[4]X[3]X[4]
  +3s[1]2s[2]2s[4]2v[0]X[0]
  \hfill\cr
\quad {}
  +(s[1]3-s[2]2s[4])s[1]s[2]2(v[6]X[1]+v[1]X[6])
  +(s[2]3-s[1]s[4]2)s[2]s[4]2(v[5]X[2]+v[2]X[5])
  +(s[4]3-s[1]2s[2])s[1]2s[4](v[4]X[3]+v[3]X[4])  \hfill\cr
\quad {}
-(4s[1]s[2]4-3s[1]2s[2]s[4]2+2s[4]5)s[1]3s[2]s[4]v[1]v[6]
  -(4s[2]s[4]4-3s[1]2s[2]2s[4]+2s[1]5)s[1]s[2]3s[4]v[2]v[5]
\hfill\cr \quad {}
-(4s[1]4s[4]-3s[1]s[2]2s[4]2+2s[2]5)s[1]s[2]s[4]3v[3]v[4]\;.
\hfill\cr}
$$
The equations can be slightly simplified by
replacing the deformation parameters  by $v_i/s_4$ (for $s_4=0$
the elliptic curve degenerates); using the equation of the modular
curve the resulting formulas are still polynomial. However such an
operation breaks the symmetry.

The base space of the versal deformation in negative degree is the cone over
an elliptic scroll of degree 7 in $\P^6$. The seven equations of the
scroll are of the following type:
$$\macaulay
s[1]s[2]s[4]v[0]2+s[1]s[4]2v[1]v[6]+s[1]2s[2]v[2]v[5]+s[4]s[2]2v[3]v[4]\;.
$$
Actually these equations do not describe the family of scrolls over the
modular curve; from each equation one can obtain three other ones
by multiplying with $s_i$, using the modular equation and dividing by
$s_{4i}$.

\subex8
To describe the deformation it is enough to give the equations
$F^0_{012}$, $F^0_{123}$, $F^0_{234}$, $F^{\ha1}_{\ha1\ha3\ha5}$
and $F^{\ha1}_{\ha3\ha5\ha7}$, which are all $F^{h}_{j-1,j,j+1}$
with $h=0$ or $h=\ha1$.
Due to the modular symmetries the following three equations suffice.
$$\macaulay
\displaylines{
\quad
s[1]2X[2]X[6]-s[2]2X[1]X[7]+s[1]s[3]X[0]2 \hfill\cr
\quad {}
  -v[0]X[0]s[2]3s[3]2+v[4](X[4]+s[2]s[4]2v[4])s[2]3s[1]2
  +(v[3]X[5]+v[5]X[3]+s[2]s[4]2v[3]v[5])s[1]3s[3]s[4] \hfill\cr
\quad {}
  +(v[7]X[1]+v[1]X[7])s[1]s[3]3s[4]
  +(v[6]X[2]+v[2]X[6]+s[2]s[4]2v[2]v[6])(s[1]2s[2]s[4]2-s[1]s[2]3s[3])
  -s[2]3(s[1]4+s[3]4)s[4]v[1]v[7] \hfill\cr
}
$$
$$\macaulay
\displaylines{
\quad
s[1]s[3]X[3]X[5]+s[1]s[3]X[1]X[7]-s[2]s[4]X[2]X[6] \hfill\cr
\quad {}
  +(v[0]X[0]+v[4]X[4])s[1]s[3]s[2](s[3]2-s[1]2)
  +(v[6]X[2]+v[2]X[6])s[2](s[3]4+s[1]4) \hfill\cr
\quad {}
  -s[2]2s[4]2s[1]s[3](s[1]2v[4]2+s[3]2v[0]2)-s[1]s[2]2s[3]s[4]4v[1]v[7]
  -s[2]3s[4]5v[2]v[6]\hfill\cr
}
$$
$$\macaulay
\displaylines{
\quad s[1]s[2]X[3]X[6]-s[2]s[3]X[2]X[7]+s[1]s[4]X[0]X[1] \hfill\cr
\quad {}
  +(v[4]X[5]+v[5]X[4])s[1]s[2]3s[4]+(v[7]X[2]+v[2]X[7])s[2]2s[3]s[4]2
  +(v[0]X[1]+v[1]X[0])(s[1]s[2]s[4]3-s[2]3s[3]s[4]) \hfill\cr
\quad {}
  -(2s[1]s[2]5s[4]2-s[2]3s[3]s[4]4+s[1]s[2]s[4]6)v[2]v[7]
  -s[2]2s[3]s[4]5v[4]v[5]\hfill\cr
}
$$
The base space (for fixed $s_\al$) is given by twelve linearly independent
quadrics, eight with odd index sum, two of which come from
the following equations with  $2h=1$:
$$\macaulay
\displaylines{
s[1]s[2]v[2]v[7]+s[3]s[2]v[3]v[6]-s[3]s[4]v[0]v[1] \cr
s[1]s[2]v[3]v[6]+s[3]s[2]v[2]v[7]-s[3]s[4]v[4]v[5] \cr
s[1]s[2]v[0]v[1]-s[3]s[2]v[4]v[5]+s[1]s[4]v[2]v[7] \cr
s[1]s[2]v[4]v[5]-s[3]s[2]v[0]v[1]+s[1]s[4]v[3]v[6] \cr
}
$$
and four even ones, which  can be written
as extension of the matrix
$$
\macaulay\pmatrix{
s[3]2-s[1]2  &s[2]s[4]     &s[1]s[3] &v[1]v[7]-v[3]v[5]\cr
s[4]2        &s[1]2+s[3]2  &s[2]2    &v[4]2-v[0]2      \cr}
$$
decribing the base curve; I have only written the $h=0$ terms.
The base space consists of five components, four planes
given by the ideals:
$$\macaulay
\displaylines{
( v[7], v[5], v[3], v[1], v[2]-v[6], v[0]-v[4])\cr
( v[7], v[5], v[3], v[1], v[2]+v[6], v[0]+v[4])\cr
( v[6], v[4], v[2], v[0], v[3]-v[7], v[1]-v[5])\cr
( v[6], v[4], v[2], v[0], v[3]+v[7], v[1]+v[5])\cr
}
$$
and one component of degree 8, which is isomorphic to
the cone over the elliptic curve. In fact, the ideal of the
component can be obtained from the ideal of the curve by the
substitution $X_i\mapsto v_i$ for even $i$ and
$X_i\mapsto v_{i+4}$ for odd $i$.

\comment
The total space of a homogeneous 1-parameter deformation
of degree $-1$ is
the cone  over a surface in $\P^8$ with the elliptic curve $E$ as
hyperplane section. In these terms the description of the
components is as follows.
Identify $\Pic^8(E)$ with ${\rm Jac}(E)\cong E$ in such a way that
$L:=\sier{}(1)$ represents the origin of the
group law.
There are 4 line bundles $M$ of degree
4 with $M^{\otimes 2}=L$. Each of these embeds $E$ in $\P^3$ and
composition with the Veronese embedding  $V\colon \P^3\to\P^9$
gives the map to $\P^7$. The image of $E$ in $\P^3$
lies on a pencil of quadrics, which is transformed by
$V$ in a pencil of hyperplane sections of $V(\P^3)$.
This construction gives  a pencil of surfaces with $E$ as
hyperplane section, and therefore a plane in the negative part of
the base  of the versal deformation. The general
quadric in the pencil is non-singular, but there are
four singular quadrics; the
resulting double cover of $\P^1$ is isomorphic to $E$.
As to the component isomorphic to the cone over $E$,  given a
point $Q\in E$  and the origin $P\in E$, I embed  $E$ in the plane
with the linear system $|2P+Q|$, then  $\P^2$ with the
triple Veronese embedding $V_3$ in $\P^9$, and finally I project the
surface $V_3(\P^2)$ from the point $V_3\circ\phi_{|2P+Q|}(3Q)$ to
$\P^8$.
\endcomment

\subex9
The equations $F^0_{ijk}$ contain (for fixed value of the modulus) three
linear independent ones; I choose $F^0_{013}$, $F^0_{043}$ and $F^0_{073}$.
As they are permuted by the action of the  modular substitution $U$ it
suffices to give the deformation of $F^0_{013}$:
$$
\macaulay
\displaylines{
\qquad
    -s[4]s[7]X[0]2+s[1]2X[3]X[6]
    -s[3]2(X[1]+s[1](s[1]2-s[4]s[7])v[1])(X[8]+s[1](s[1]2-s[4]s[7])v[8])
\hfill\cr\qquad{}
   -s[7]s[4](s[4]3+3s[1]3+2s[7]3)v[0]X[0]
    +s[1]2(s[7]3+2s[4]3)(v[6]X[3]+v[3]X[6])
\hfill\cr\qquad{}
    +s[3]2s[1](s[1]s[7]-s[4]2)(v[7]X[2]+v[2]X[7])
    +s[3]2s[1](s[4]s[1]-s[7]2)(v[5]X[4]+v[4]X[5])
\hfill\cr\qquad{}
        -s[4]s[7](s[4]3+s[1]3+s[7]3)(2s[1]3+s[7]3)v[0]2
    +3s[1]2s[4]3(2s[7]3+s[4]3)v[3]v[6]
\hfill\cr\qquad{}
   +s[3]2s[1](s[1]s[7]-s[4]2)(s[1]3+2s[4]3-s[1]s[4]s[7]+s[7]3)v[2]v[7]
   +s[3]2s[7](s[4]5-3s[1]3s[4]2+8s[1]2s[4]s[7]2+4s[4]2s[7]3-s[1]s[7]4)v[4]v[5]
\cr}
$$
The base space is given by nine times two equations:
$$
\macaulay
\displaylines{
 s[3](v[0]2-v[3]v[6])-s[4]v[1]v[8]-s[1]v[7]v[2]-s[7]v[4]v[5] \cr
 s[1]v[1]v[8]+s[7]v[7]v[2]+s[4]v[4]v[5]\;.\cr
}
$$
In the actual computation these equations come out multiplied with some
factors in the modular parameters $s_i$, which vanish at the cusps. As the
elliptic curve is supposed to be smooth there is no harm in dividing by those
factors. One of them is the Hessian $
\macaulay H=s[1]3+s[4]3+s[7]3-3s[1]s[4]s[7]$.
The base space has nine irreducible components, three of which are given by
the ideal
$$\macaulay
(v[1],v[4],v[7],v[2],v[5],v[8], v[0]2-v[3]v[6],v[3]2-v[0]v[6], v[6]2-v[0]v[3])
$$
and the other six by cyclic permutation. Furthermore there is an embedded
component. The equations for the reduction of the base space include $v_iv_j$
if $i\not\equiv j \pmod3$, which simplifies the equation $F^0_{013}$
considerably.

The total space over each component is the cone over the triple Veronese
embedding of $\P^2$. I give the formulas for the case $v_0=v_3=v_6$.
The projection
$$
\macaulay
\displaylines{
\qquad
Y[0] = s[3](s[7]-s[4])(X[1]+X[8])+s[3](s[4]-s[1])(X[2]+X[7])
     +s[3](s[1]-s[7])(X[4]+X[5]) \hfill\cr \hfill{}
    -(s[1]s[4]+s[4]s[7]+s[1]s[7])(X[3]+X[6])-(s[1]2+s[4]2+s[7]2)X[0]
\qquad\cr\qquad
Y[3] = s[3](s[4]-s[1])(X[1]+X[5])+s[3](s[1]-s[7])(X[8]+X[7])
    +s[3](s[7]-s[4])(X[4]+X[2])  \hfill\cr \hfill{}
    -(s[1]s[4]+s[4]s[7]+s[1]s[7])(X[0]+X[6])-(s[1]2+s[4]2+s[7]2)X[3]
\qquad\cr\qquad
Y[6] = s[3](s[4]-s[1])(X[1]+X[2])+s[3](s[4]-s[1])(X[4]+X[8])
    +s[3](s[7]-s[4])(X[5]+X[7])  \hfill\cr \hfill{}
    -(s[1]s[4]+s[4]s[7]+s[1]s[7])(X[3]+X[0])-(s[1]2+s[4]2+s[7]2)X[6]\;,
\qquad\cr}
$$
from the space spanned by six suitable chosen points
maps the elliptic curve $E$ onto the plane cubic
$$ \macaulay
f=s[3]3(Y[0]3+Y[3]3+Y[6]3)+(s[1]3+s[4]3+s[7]3+6s[1]s[4]s[7])Y[0]Y[3]Y[6] \;.
$$
The functions $Y_l$ can be obtained from the sum $\sum X_iX_jX_k$,
where the sum ranges over  $i,j,k$ in different residue classes modulo 3
with $i+j+k\equiv l\pmod 9$,
by specialising the first two factors of each product.
The equation $f$ together with the following nine polynomials
$$
\macaulay
\eqalign{
X[0] {}&= -(s[3]3Y[0]2+3s[1]s[4]s[7]Y[3]Y[6])Y[0]{} \cr
X[3] {}&= -(s[3]3Y[3]2+3s[1]s[4]s[7]Y[0]Y[6])Y[3]{} \cr
X[6] {}&= -(s[3]3Y[6]2+3s[1]s[4]s[7]Y[3]Y[0])Y[6]{} \cr
X[1] {}&= s[3]2(s[4]Y[0]2Y[3]+s[1]Y[3]2Y[6]+s[7]Y[6]2Y[0]) \cr
X[4] {}&= s[3]2(s[7]Y[0]2Y[3]+s[4]Y[3]2Y[6]+s[1]Y[6]2Y[0]) \cr
X[7] {}&= s[3]2(s[1]Y[0]2Y[3]+s[7]Y[3]2Y[6]+s[4]Y[6]2Y[0]) \cr
X[2] {}&= s[3]2(s[1]Y[0]2Y[6]+s[4]Y[3]2Y[0]+s[7]Y[6]2Y[3]) \cr
X[5] {}&= s[3]2(s[7]Y[0]2Y[6]+s[1]Y[3]2Y[0]+s[4]Y[6]2Y[3]) \cr
X[8] {}&= s[3]2(s[4]Y[0]2Y[6]+s[7]Y[3]2Y[0]+s[1]Y[6]2Y[3]) \cr
}
$$
forms a basis of the linear system of plane cubics, except when
$\macaulay H=s[1]3+s[4]3+s[7]3-3s[1]s[4]s[7] =0$ or $s_3=0$, which are
precisely the equations of the cusps of the modular curve.
Upon setting $v_0=f/H$ the image of $\P^2$ exactly satisfies the equations
above.

The group of points of order nine of the elliptic curve  acts on the set of
components with a points of order three giving rise to an  automorphism of
each of the components of the total space.

\biblio

\item \rfarn
  J\"urgen Arndt,
  {\sl  Verselle Deformationen zyklischer
  Quotientensingularit\"aten\/}.
  Diss. Hamburg 1988.

\item\rfbh
  Wolf Barth and Klaus Hulek,
  {\sl Projective Models of Shioda Modular Surfaces\/},
  manu\-scrip\-ta math. {\bf 50} (1985), 73--112.

\item \rfbar
  R. Bartsch,
  {\sl Meromorphe Funktionen auf der universellen elliptischen Kurve
  mit Niveau $N$-Struktur\/}.
  Diss. Hamburg 1985.

\item \rfbehe
  Kurt Behnke,
  {\sl Infinitesimal Deformations of Cusp
  Singularities\/}.
  Math. Ann. {\bf 265} (1983), 407--422.

\item \rfbeht
  Kurt Behnke,
  {\sl On the Module of Zariski Differentials and Infinitesimal
  Deformations of Cusp Singularities\/}.
  Math. Ann. {\bf 271} (1985), 133--142.

\item \rfbi
  L. Bianchi,
  {\sl \"Uber die Normalformen dritter und f\"unfter Stufe des
  elliptischen Integrals erster Gattung\/}.
  Math. Ann. {\bf 17} (1880), 234--262.

\item \rfcoh
  Harvey Cohn,
  {\sl Some Explicit Resolutions of Modular Cusp
  Singularities\/}.
  Math. Ann. {\bf 198} (1972), 123--130.

\item \rfdr
  P. Deligne et M. Rappaport,
  {\sl Les sch\'emas de modules de courbes elliptiques\/}.
  In: Modular Functions of One Variable II,
  Berlin etc., Springer 1973 (Lect. Notes in Math.; 349),
  pp. 143--316.

\item\rfgeer
  Gerard van der Geer,
  {\sl Hilbert Modular Surfaces\/}.
  Berlin etc., Springer 1988.

\item\rfhal
  G.-H. Halphen,
  {\sl Trait\'e des fonctions elliptiques et de leurs applications\/},
  Paris, Gauthier-Villars 1886.

\item\rfhir
  F. Hirzebruch,
  {\sl Hilbert modular surfaces\/}.
  L'Ens. Math {\bf 71} (1973) 183--281.

\item \rfhul
  Klaus Hulek,
  {\sl Projective geometry of elliptic curves\/}.
  Ast\'erisque {\bf 137} (1986).

\item \rfhur
  Adolf Hurwitz,
  {\sl \"Uber endliche Gruppen linearer Substitutionen,
  welche in der Theorie der elliptischen Transcendenten auftreten\/}.
  Math. Ann. {\bf 27} (1886), 183--233.

\item \rfhurf
  Adolf Hurwitz,
  {\sl \"Uber die Weierstrass'sche $\sigma$-Funktion\/}.
  In: Festschrift f\"ur H.A. Schwarz, Berlin 1914, pp. 133-141.
  Also in: Mathematische Werke Bd. I, pp. 722-730.

\item \rfkm
  Nicholas M. Katz and Barry Mazur,
  {\sl Arithmetic moduli of elliptic curves\/}.
  Princeton 1985 (Ann. of Math. Studies; 108).

\item \rfklel
  Felix Klein,
  {\sl \"Uber die elliptischen Normalkurven der $n$-ten
  Ordnung\/}. Nr. XC in: Gesammelte Mathematische
  Abhandlungen Bd. III, Berlin, Springer 1923,
  pp. 198--254.

\item \rfkf
  Felix Klein,
  {\sl Vorlesungen \"uber die Theorie der elliptischen
  Modulfunctionen\/}, ausgearbeitet und vervollst\"andigt
  von Dr. Robert Fricke. Zweiter Band,
  Leipzig, B. G. Teubner 1892.

\item \rfmer
  J.-Y. M\'erindol,
  {\sl Les singularit\'es simples elliptiques, leur d\'eformations,
  les surfaces de Del Pezzo et les transformations  quadratiques}.
  Ann. scient. \'Ec. Norm. Sup. {\bf 15} (1982), 17--44.

\item \rfmf
  D. Mumford and J. Fogarty,
  {\sl Geometric Invariant Theory\/},
  Second Enlarged Edition.
  Berlin etc., Springer 1982.

\item \rfpe
  Ruud Pellikaan, {\sl Personal communication}.

\item \rfpr
  J\"urgen Pesselhoy and Oswald Riemenschneider,
  {\sl Projective resolutions of Hodge algebras: some
  examples\/}.
  Proc. Symp. Pure Math. {\bf 40}, Part 2 (1983), 305--317.

\item \rfrie
  Oswald Riemenschneider,
  {\sl Deformationen von Quotientensingularit\" aten {\rm(}nach
  zykli\-schen Gruppen\/}).
  Math. Ann. {\bf209} (1974), 211--248.

\item  \rfrie
  Oswald Riemenschneider,
  {\sl Zweidimensionale  Quotientensingularit\"aten\/{\rm :}
  Gleichungen\break und Syzygien\/}.
  Arch. Math. {\bf 37} (1981), 406--417.

\item \rfvel
  Jacques V\'elu,
  {\sl Courbes elliptiques munies d'un sous-groupe $\Z/n\Z \times
  \mu_n$\/}.
  Bull. Soc. Math. France, M\'emoire {\bf 57} (1978), 5--152.

\item \rfwei
  Karl Weierstra\ss,
  {\sl Zur Theorie der Jacobi'schen Functionen von mehreren
  Ver\"ander\-lichen}.
  Sitzungsber. K\"onigl. Preuss. Akad. der Wiss. (1882) =
  Werke, Bd. III, pp. 155-159.

\vfill
\parindent=0pt
{\obeylines
Address of the author:
Matematiska Institutionen, Chalmers Tekniska H\"ogskola
S 41296 G\"oteborg, Sweden
E-mail:  stevens@math.chalmers.se     }

\bye